\documentclass[conference]{IEEEtran}
\usepackage{cite}
\usepackage[]{inputenc}
\usepackage{textcomp}
\def\BibTeX{{\rm B\kern-.05em{\sc i\kern-.025em b}\kern-.08em
    T\kern-.1667em\lower.7ex\hbox{E}\kern-.125emX}}
    
\makeatletter

\title{\huge{User-Device Authentication in Mobile Banking \\using APHEN for Paratuck2 Tensor Decomposition}}

\usepackage{amsmath}
\usepackage{amsfonts}
\usepackage{epsfig}
\usepackage[mathscr]{euscript}
\usepackage{booktabs} 
\usepackage{stmaryrd} 
\usepackage{hyperref} 
\usepackage{graphics} 
\usepackage{bbm} 
\usepackage{blindtext}
\usepackage{subcaption}
\usepackage{tikz}
\usepackage{standalone}
\usetikzlibrary{arrows, patterns}
\usetikzlibrary{decorations.pathreplacing,angles,quotes}
\usepackage{framed, standalone}
\usepackage{capt-of}
\usepackage{balance}
\usepackage{makecell} 
\usepackage[linesnumbered,lined,boxed,commentsnumbered,ruled,vlined]{algorithm2e}
\usepackage{setspace}

\begin{document}



\author{\IEEEauthorblockN{Jeremy Charlier}
\IEEEauthorblockA{
University of Luxembourg\\
Luxembourg \\
jeremy.charlier@uni.lu}
\and
\IEEEauthorblockN{Eric Falk}
\IEEEauthorblockA{
University of Luxembourg\\
Luxembourg \\
eric.falk@uni.lu}
\and
\IEEEauthorblockN{Radu State}
\IEEEauthorblockA{
University of Luxembourg\\
Luxembourg \\
radu.state@uni.lu}
\and
\IEEEauthorblockN{Jean Hilger}
\IEEEauthorblockA{BCEE\\
Luxembourg \\
j.hilger@bcee.lu}
}

\maketitle

\begin{abstract}
The new financial European regulations such as PSD2 are changing the retail banking services. Noticeably, the monitoring of the personal expenses is now opened to other institutions than retail banks. Nonetheless, the retail banks are looking to leverage the user-device authentication on the mobile banking applications to enhance the personal financial advertisement. To address the profiling of the authentication, we rely on tensor decomposition, a higher dimensional analogue of matrix decomposition. We use Paratuck2, which expresses a tensor as a multiplication of matrices and diagonal tensors, because of the imbalance between the number of users and devices. We highlight why Paratuck2 is more appropriate in this case than the popular CP tensor decomposition, which decomposes a tensor as a sum of rank-one tensors. However, the computation of Paratuck2 is computational intensive. We propose a new APproximate HEssian-based Newton resolution algorithm, APHEN, capable of solving Paratuck2 more accurately and faster than the other popular approaches based on alternating least square or gradient descent. The results of Paratuck2 are used for the predictions of users' authentication with neural networks. We apply our method for the concrete case of  targeting clients  for financial advertising campaigns based on the authentication events generated  by mobile banking applications.
%
%
%
%
%
\end{abstract}

%
%
%
\section{Introduction}
Endorsed by the  European objectives to promote the financial exchanges between the Euro members, a new financial regulation for Personal Finance Management (PFM) has entered into force in 2018. PFM is the monitoring of the revenues and the expenses of a bank account. PFM is achieved with the use of a Personal Finance Application (PFA), otherwise known as mobile banking application. The revised Payment Service Directive, PSD2, allows every person having a bank account to use a PFA from a third party provider to manage its personal finance, and thus transform the banks into simple vaults. 
Nevertheless, the retail banks now have the opportunity to leverage the user-device authentication on their own mobile banking application. Through the regular authentication, the bank create a financial profile awareness for every clients. They measure the frequency of the connections per day, the time of the connections and the type of the mobile devices used to authenticate such as a smartphone or a tablet. The more frequently a client is authenticating to its mobile banking application, the more likely he will have a high interest for finance, and therefore, the more likely he will be interested by financial recommendation. This client will be contacted in priority to advert financial products for wealth optimization. However, it is very common for the same person to use several mobile devices for the same application. Therefore, the devices generate dozens of authentication per day and create a strong imbalance between the number of users and the number of devices. Consequently, the modeling of the user-device authentication is multidimensional and complex. To answer this challenge, we rely on tensor decompositions, a higher order analogue of matrix decompositions, since they have proven to be powerful to model multidimensional interactions. Particularly, we address the modeling of the imbalanced user-device authentication with the Paratuck2 tensor decomposition, which decomposes a tensor as a multiplication of matrices and sparse tensors. We summarize the three main contributions of the paper as follows:

\begin{itemize}
\setlength{\itemsep}{0pt}%
\setlength{\parskip}{0pt}%
\vspace{-0.1cm}
\item We have designed an innovative AProximate HEssian Newton minimization algorithm, APHEN, applied to a tensor decomposition, Paratuck2. The algorithm is able to reduce at a minimum the numerical errors while delivering state of the art performance. Additionally, APHEN does not require the full knowledge of the Hessian matrix to achieve its superior convergence.

\item We highlight the superior capabilities of APHEN and Hessian-based resolution algorithms for complex tensor decomposition such as Paratuck2 in comparison to other popular resolution schemes, mainly exclusively studied for the popular CP decomposition. 


\item  Additionally, we have applied the unique advantages of Paratuck2 for interactions modeling with imbalanced data. We justify with numerical simulation the use of Paratuck2 instead of the more popular CP tensor decomposition, which decomposes a tensor as a sum of rank-one tensors. In our application, one user can use several computers, and thus, we have considered the imbalance between the number of users and computers.

\item Finally, we have developed an approach with the Paratuck2 tensor decomposition and neural networks for authentication monitoring. Based on Paratuck2, the neurons predict the users' authentication to estimate the future financial awareness of the clients. Therefore, the banks can better advertise their products by contacting the clients which are the more likely to be interested.
\end{itemize}

The remaining of this paper is organized as follows. Section 2 surveys the latest research publications related to user-device authentication and tensor decompositions. Section 3 describes the fundamentals of tensor decompositions, the APHEN algorithm and the other popular resolution algorithms. Additionally, it introduces a basic knowledge of neural networks and machine learning predictions. Section 4 illustrates the convergence speeds of APHEN in comparison to the other popular schemes. Then, Paratuck2 and neural networks are used to predict the imbalanced users' authentication with the aim of improving the subscription rates to financial products during the banks' advertising campaigns. Finally, we conclude the paper and we highlight pointers to future works in the last section.

%
%
%
\section{Related Work}
%
%
%
%
\textbf{Literature on User-Device Authentication}
The user-device authentication has significantly evolved for the past few years thanks to the new technologies. A reliable user-device authentication was proposed in \cite{skinner2016cyber}, based on a graphical user friendly authentication. In \cite{adeka2017africa}, the use of Location Based Authentication (LBA) was studied. The development of recent embedded systems within smart-devices leads to new authentication processes which were considered as a pure fiction only few years ago. In \cite{li2016face}, the usage of the embedded camera of smart-devices for authentication by face recognition was assessed. The face image taken by the camera of the mobile device was sent to a background server to perform the calculation which reverts then to the phone. In a similar approach, the use of iris recognition was proposed in \cite{chhabra2016low}. However, the authors showed this kind of authentication was not the preferred choice of the end user. Additionally, the sensors embedded into smart-devices allow other type of biometric authentication. In \cite{spolaor2016biometric}, the different biometric authentication that could be used with smart-devices were presented, such as the pulse-response, the fingerprint or even the ear shape. Although biometric or LBA solutions might offer a higher level of security for authentication, their extension toward a large scale usage is complex. In \cite{theofanos2016secure}, the authors developed the idea that public-key infrastructure based systems, such as strong passwords in combination with physical tokens, for example, a cell phone, would be more likely to be used and largely deployed. Nonetheless, it is worth mentioning that the most common procedure for mobile devices authentication is still a code of four or six digits \cite{bultel2018security}. \\




\textbf{Literature on Tensor Decomposition}
The interactions modeling of the user-device authentication is multidimensional and complex, both specificities of tensor analysis and tensor decompositions. Tensor decompositions are an extension to higher dimensions of the two dimensional matrix decompositions \cite{harshman1970foundations,carroll1970analysis}. It is explained by the evolution towards more extensive analysis in the presence of an increasing number of features within the datasets. As a result, different tensor decompositions, or tensor factorizations, exist with different resolution algorithms for different types of applications \cite{kolda2009tensor, acar2011all}. Meanwhile, the scope of the tensor application have skyrocketed. In \cite{da2016tensor}, the CP tensor decomposition was used for data mining and signal processing. A low-rank approximation  was developed for fast computation that could be used on large datasets. The CP tensor decomposition was also used in location-based social networks for identification profiling \cite{papalexakis2014spotting}. In the experiments, a certain number of anomalies were identified in the check-in behavior of the users. Following the trend of social network studies, the algorithm Tensorcat was specifically designed to study interactions on social networks \cite{de2017tensorcast}. Different sources were incorporated in coupled tensors for time-evolving networks such as Twitter to predict the evolution of the social network activity. Tensor predictive analytics have also been addressed in \cite{takeuchi2017autoregressive}. A tensor factorization method is described for spatial and temporal autocorrelations. The analysis and the predictions are demonstrated on traffic transporting data. Similarly, the Rescal tensor decomposition has been used in \cite{wang2016learning} for the review of spam detection. The approach highlighted the interactions between the reviewers and the products, and it led to a better accuracy of spam detection when compared to other methods. With similar objectives, a specific algorithm was developed in \cite{qiao2016signal} for heterogeneous data relying on the Higher-Order Singular Value Decomposition (HOSVD) to describe the frequencies of various signals. 
To answer to the real-world problems with the velocity of streaming data, multi-aspect streaming tensor completion is underlined in \cite{song2017multi}. Albeit the approach allows to build dynamic tensors, it relies on the CP tensor decomposition which lacks the linear independence of the latent variables in each order \cite{acar2011scalable}.\\

In this paper, we extend the state of the art of the tensor numerical resolution introduced for the CP decomposition in 
\cite{acar2011all, paatero1997weighted} and \cite{tomasi2006comparison}
by proposing APHEN, an APproximate HEssian Newton resolution that does not require the complete knowledge of the Hessian matrix, and therefore removing the limitations of the computational cost of the Hessian matrix. APHEN is capable of minimizing at a minimum the numerical error at convergence while having similar or faster computation time than other popular resolution schemes. Additionally, we highlight experimentally the limitations of CP \cite{acar2011scalable} and we propose the use of Paratuck2 for imbalanced data. Finally, we rely on neural networks for the predictions of the users' authentication for personalized financial recommendation occurring during the advertising campaigns relying on mobile banking application connections. 


%
%
%
\section{Model Description}
In this section, we describe both the CP and the Paratuck2 tensor decompositions initially introduced in \cite{harshman1970foundations,carroll1970analysis} and \cite{harshman1996uniqueness}. Subsequently, we describe the error minimization algorithm APHEN, which is at the core of our contribution. 
Finally, we briefly describe neural networks applied to Paratuck2 for the aim of the latent users' authentication predictions.

\subsection{CP and Paratuck2 Tensor Decompositions}\label{subsec:para2}
\noindent \textbf{Notation} The terminology hereinafter follows the one described by Kolda and Bader in \cite{kolda2009tensor} and commonly used. Scalars are denoted by lower case letters, \textit{a}. Vectors and matrices are described by boldface lowercase letters and boldface capital letters, respectively \textbf{a} and \textbf{A}. High order tensors are represented using Euler script notation, $\mathscr{X}$. The transpose matrix of $A\in \mathbb{R}^{I\times J}$ is denoted by $A^T$. The inverse of a matrix $A\in \mathbb{R}^{I\times I}$ is denoted by $A^{-1}$.

\noindent \textbf{Algebra Operations} 
The outer product between two vectors, $\textbf{u}\in\mathbb{R}^I$ and $\textbf{v}\in\mathbb{R}^J$ is denoted by the symbol $\circ$.
\begin{equation}
\textbf{u} \circ \textbf{v} = 
\begin{bmatrix}
u_1v_1 & u_1v_2 & \cdots & u_1v_J\\ 
u_2v_1 & u_2v_2 & \cdots & u_2v_J\\ 
\vdots & \vdots & \vdots & \vdots \\
u_Iv_1 & u_Iv_2 & \cdots & u_Iv_J 
\end{bmatrix}
= \textbf{u}_i\textbf{v}_j
\end{equation}

\noindent The Kronecker product between two matrices \textbf{A}$\in\mathbb{R}^{I\times J}$ and \textbf{B}$\in\mathbb{R}^{K\times L}$, denoted by \textbf{A}$\otimes$\textbf{B}, results in a matrix \textbf{C}$\in\mathbb{R}^{IK\times KL}$.
\begin{equation} \label{eq::kron}
\textbf{C}=\textbf{A}\otimes\textbf{B}=
\begin{bmatrix}
 a_{11}\textbf{B}& a_{12}\textbf{B} & \cdots & a_{1J}\textbf{B}\\ 
 a_{21}\textbf{B}& a_{22}\textbf{B} & \cdots & a_{2J}\textbf{B}\\
 \vdots & \vdots & \ddots & \vdots \\ 
 a_{I1}\textbf{B}& a_{I2}\textbf{B} & \cdots & a_{IJ}\textbf{B}
\end{bmatrix}
\end{equation}

\noindent The Khatri-Rao product between two matrices \textbf{A}$\in\mathbb{R}^{I\times K}$ and \textbf{B}$\in\mathbb{R}^{J\times K}$, denoted by \textbf{A}$\odot$\textbf{B}, results in a matrix \textbf{C} of size $\mathbb{R}^{IJ\times K}$. It is the column-wise Kronecker product.
\begin{equation} \label{eq::kr}
\textbf{C}=\textbf{A}\odot\textbf{B}=
[\textbf{a}_1\otimes \textbf{b}_1 \quad \textbf{a}_2\otimes \textbf{b}_2 \quad \cdots \quad \textbf{a}_K\otimes \textbf{b}_K]
\end{equation}

\noindent \textbf{Tensor Definition} $\mathscr{X}$ is called a \textit{n}-way tensor if  $\mathscr{X}$ is a \textit{n}-th multidimensional array. It is expressed as $\mathscr{X}\in \mathbb{R}^{I_1 \times I_2 \times I_3 \times ... \times I_n}$.

\noindent \textbf{Tensor Operations} The square root of the sum of all tensor entries squared of the tensor $\mathscr{X}$ defines its norm.
\begin{equation} \label{eq::norm}
||\mathscr{X}||=\sqrt{\sum_{j=1}^{I_1}\sum_{j=2}^{I_2}...\sum_{j=n}^{I_n}x_{j_1, j_2, ..., j_n}^2}
\end{equation}
The rank-\textit{R} of a tensor $\mathscr{X}\in\mathbb{R}^{I_1\times I_2\times ...\times I_N}$ is the number of linear components that could fit $\mathscr{X}$ exactly.
\begin{equation} \label{eq::rank}
\mathscr{X}=\sum_{r=1}^R \textbf{a}_r^{(1)} \circ \textbf{a}_r^{(2)} \circ ... \circ \textbf{a}_r^{(N)}
\end{equation}

\noindent \textbf{Definition of the CP Decomposition} We motivate the use of Paratuck2 over CP because of the imbalance in our dataset. Effectively, CP lacks the linear independence of the factors in each order \cite{acar2011scalable}. The CP decomposition has been introduced in \cite{harshman1970foundations, carroll1970analysis}. The tensor $\mathscr{X}\in\mathbb{R}^{I\times I\times K}$ is defined as a sum of rank-one tensor. The number of rank-one tensors is determined by the rank, denoted by $R$, of the tensor $\mathscr{X}$. The CP decomposition is expressed as 
\begin{equation}
\label{eq::cpequation}
\mathscr{X} = \sum_{r=1}^{R} \textbf{a}_r^{(1)} \circ \textbf{a}_r^{(2)} \circ \textbf{a}_r^{(3)} \circ... \circ \textbf{a}_r^{(N)}
\end{equation}
where $\textbf{a}_r^{(1)}, \textbf{a}_r^{(2)}, \textbf{a}_r^{(3)}, ..., \textbf{a}_r^{(N)}$ are vectors of size $\mathbb{R}^{I_1}, \mathbb{R}^{I_2}, \mathbb{R}^{I_3}, ..., \mathbb{R}^{I_N}$. Each vector $\textbf{a}_r^{(n)}$ with $n\in \left\lbrace 1, 2, ..., N \right\rbrace$ and $r \in \left\lbrace 1, ..., R \right\rbrace$ refers to one dimension and one rank of the tensor $\mathscr{X}$.

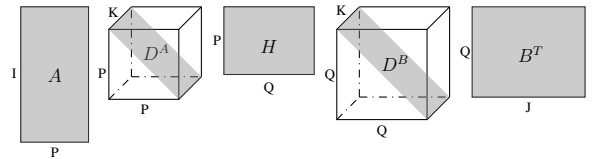
\begin{figure}[b]
\begin{center}
\begin{tikzpicture}[scale=0.6, transform shape]
\draw  (-4.5,2) rectangle (-3,-1);
\fill [gray, opacity=.4]  (-4.5,2) rectangle (-3,-1);
\draw (-3.75,0.75) node[below]{\large{$A$}};

\draw (2.5,-0.5) rectangle (4.5,1.5) node (v1) {};
\draw (4.5,1.5) -- (5,2);
\draw (2.5,1.5) -- (3,2);
\draw (4.5,-0.5) -- (5,0);
\draw (3,2) -- (5,2);
\draw (5,0) -- (5,2);
\draw [dashdotted] (3,0) -- (5,0);
\draw [dashdotted] (3,0) -- (3,2);
\draw [dashdotted] (2.5,-0.5) -- (3,0);
\fill [gray, opacity=.4] (4.5,-0.5) -- (5,0) -- (3,2) -- (2.5,1.5)  -- cycle;
\draw (-1.5,1.3) node[below]{\large{$D^A$}};

\draw (0,2) rectangle (2,0.5);
\fill [gray, opacity=.4] (0,2) rectangle (2,0.5);
\draw (1,1.4) node[below]{\large{$H$}};

\draw (-2.55,-0.05) rectangle (-1,1.5) node (v1) {};
\draw (-1,1.5) -- (-0.5,2);
\draw (-2.55,1.5) -- (-2.05,2);
\draw (-1,-0.05) -- (-0.5,0.45);
\draw (-2.05,2) -- (-0.5,2);
\draw (-0.5,0.45) -- (-0.5,2);
\draw [dashdotted] (-2.05,0.45) -- (-0.5,0.45);
\draw [dashdotted] (-2.05,0.45) -- (-2.05,2);
\draw [dashdotted] (-2.55,-0.05) -- (-2.05,0.45);
\fill [gray, opacity=.4] (-1,-0.05) -- (-0.5,0.45) -- (-2.05,2.0) -- (-2.55,1.5)  -- cycle;
\draw (3.8,1.05) node[below]{\large{$D^B$}};

\draw  (5.5,2) rectangle (8,0);
\fill [gray, opacity=.4] (5.5,2) rectangle (8,0);
\draw (6.8,1.25) node[below]{\large{$B^T$}};

\draw (-4.65,0.75) node[below]{\small{I}};
\draw (-3.75,-1) node[below]{\small{P}};
\draw (-1.75,-0.05) node[below]{\small{P}};
\draw (-2.7,0.75) node[below]{\small{P}};
\draw (-2.45,2.1) node[below]{\small{K}};
\draw (2.6,2.15) node[below]{\small{K}};
\draw (-0.15,1.5) node[below]{\small{P}};
\draw (1,0.45) node[below]{\small{Q}};
\draw (2.35,0.75) node[below]{\small{Q}};
\draw (3.5,-0.5) node[below]{\small{Q}};
\draw (5.35,1.25) node[below]{\small{Q}};
\draw (6.75,0) node[below]{\small{J}};
\end{tikzpicture}
\caption{Paratuck2 decomposition of a three-way tensor with dimension notations}
\label{fig::PARATUCK2}
\end{center}
\end{figure}

\noindent \textbf{Definition of the Paratuck2 Decomposition} Paratuck2 has been introduced by Harshman and Lundy in \cite{harshman1996uniqueness}. The tensor $\mathscr{X}\in\mathbb{R}^{I\times J\times K}$ is described  as a product of matrices and tensors
\begin{equation} \label{eq::paratuck2}
\textbf{X}_k = \textbf{A}\textbf{D}^A_k\textbf{H}\textbf{D}^B_k\textbf{B}^T \quad \textrm{with} \quad k=\left\lbrace 1, ..., K \right\rbrace
\end{equation}
where $\textbf{A}$, $\textbf{H}$ and $\textbf{B}$ are matrices of size $\mathbb{R}^{I\times P}$, $\mathbb{R}^{P \times Q}$ and $\mathbb{R}^{J\times Q}$.
The matrices $\textbf{D}^A_k\in \mathbb{R}^{P\times P}$ and $\textbf{D}^B_k\in \mathbb{R}^{Q\times Q}\:\forall k\in\left\lbrace 1, ...,K \right\rbrace$ are the slices of the tensors $\mathscr{D}^A\in \mathbb{R}^{P\times P \times K}$ and $\mathscr{D}^B\in \mathbb{R}^{Q\times Q \times K}$. The latent factors $P$ and $Q$ are related to the rank of each object set as illustrated in figure \ref{fig::PARATUCK2}.
More precisely, the columns of the matrices $\textbf{A}$ and $\textbf{B}$ represent the latent factors $P$ and $Q$.
The matrix $\textbf{H}$ underlines the asymmetry between the $P$ latent factors and the $Q$ latent factors. The tensors $\mathscr{D}^A$ and $\mathscr{D}^B$ measures the evolution of the  latent factors regarding the third dimension.

\subsection{APHEN and Approximate Derivatives}
The Alternating Least Square (ALS) method is the most commonly used method for tensor resolution as initially described in \cite{harshman1970foundations, carroll1970analysis}.
It has been applied by Bro in \cite{bro1998multi} to Paratuck2. Nonetheless, with larger data sets, the convergence performance of the ALS method decreases. To overcome this, algorithms facing gradient resolution for tensors have emerged \cite{paatero1997weighted, acar2011scalable, acar2011all}. However, it is well known that gradient descent schemes are very sensitive to the initial guess and the local minima. Therefore, we propose ApHeN, a resolution scheme that relies on the Newton conjugate gradient but does not require the knowledge of the complete Hessian matrix. Last but not least, the algorithm is applied to the Paratuck2 tensor decomposition illustrated in figure \ref{fig::PARATUCK2}.

The objective minimization function is denoted by $f$.
\begin{equation} \label{eq::errorminimization}
f(\textbf{x}) = \min_{\mathscr{\hat{X}}} ||\mathscr{X}-\mathscr{\hat{X}}||
\end{equation}
The tensor $\mathscr{\hat{X}}$ is the approximate tensor of $\mathscr{X}$ built from the decomposition with matrices $\textbf{A}$, $\textbf{H}$ and $\textbf{B}$ initially randomized.
The diagonal entries of the tensors $\mathscr{D}^A$ and $\mathscr{D}^B$ are set to 1 at the beginning of the minimization process.

\noindent The vector $\textbf{x}$ is a flattened vector containing all the entries of involved in the decomposition scheme to build $\mathscr{\hat{X}}$.
\begin{equation} \label{eq::vectorize_x}
\textbf{x} =
\begin{bmatrix}
a_{11} \; a_{12} \cdots 
a_{IP} \; d^a_{111} \cdots 
d^a_{PPK} \; d^b_{111} \cdots 
b_{JQ} 
\end{bmatrix}^T
\end{equation}

\noindent Using the notation from equation \ref{eq::vectorize_x}, we can derive the gradient and the Hessian matrix related to the CP decomposition.

The gradient, denoted by $\nabla f$, is a vector containing all the first derivatives of the function $f$ with respect to $\textbf{x}$.
\begin{equation} \label{eq::gradient}
\nabla f = 
\textrm{grad}_x f = 
\begin{bmatrix}
\dfrac{\partial f}{\partial x_1}, 
\dfrac{\partial f}{\partial x_2}, 
\cdots, 
\dfrac{\partial f}{\partial x_n}
\end{bmatrix}
\end{equation}
The Hessian matrix, $\textbf{Hes}$, is the matrix containing the second derivatives of the function $f$ with respect to $\textbf{x}$.
\begin{equation} \label{eq::Hessian_matrix}
\textbf{Hes} = \begin{pmatrix}
\dfrac{\partial f}{\partial x_1^2} & 
\dfrac{\partial f}{\partial x_1 \partial x_2} & 
\cdots & 
\dfrac{\partial f}{\partial x_1 \partial x_n} \\[1em]
\dfrac{\partial f}{\partial x_2 \partial x_1} & 
\dfrac{\partial f}{\partial x_2^2} & 
\cdots & 
\dfrac{\partial f}{\partial x_2 \partial x_n} \\
\vdots & \vdots & \ddots & \vdots \\
\dfrac{\partial f}{\partial x_n \partial x_1} &
\dfrac{\partial f}{\partial x_n \partial x_2} & 
\cdots & 
\dfrac{\partial f}{\partial x_n^2}
\end{pmatrix}
\end{equation}
The Newton conjugate gradient algorithm minimizes the function $f$ according to the below equation.
\begin{equation}
\begin{aligned}
f(\textbf{x}) = 
f(\textbf{x}_0)
& + \nabla f(\textbf{x}_0) (\textbf{x}-\textbf{x}_0) \\
& +\dfrac{1}{2}(\textbf{x}-\textbf{x}_0)^T 
\textbf{Hes}(\textbf{x}_0)(\textbf{x}-\textbf{x}_0)
\end{aligned}
\end{equation}

The variable $\textbf{x}_0$ is the initial guess, $\nabla f$ the gradient of $f$ and $\textbf{Hes}$ the Hessian matrix of $f$.
If the Hessian matrix is positive definite then the local minimum of the function is determined by setting the gradient of the quadratic form to zero. 
\begin{equation} \label{eq::update_x}
\textbf{x}_{opt} = \textbf{x}_0 - \textbf{Hes}^{-1} \nabla f
\end{equation}

Since the gradient and the Hessian matrix are computed with finite differences, the only prerequisite for Paratuck2 tensor decomposition is the factorization equation (\ref{eq::paratuck2}). Thus, the method can be transposed to other decompositions, such as CP, by merely changing the tensor decomposition equation. The approximate gradient is based on a fourth order formula (\ref{eq::fourth_order_diff}) to ensure reliable approximation \cite{schittkowski2001nlpqlp}.
\begin{equation} \label{eq::fourth_order_diff}
\begin{split}
\dfrac{\partial}{\partial x_i} f(\textbf{x}) 
\approx 
\dfrac{1}{4!\eta} & \big( \; 2f(\textbf{x}-2\eta e_i) 
-16 f(\textbf{x} - \eta e_i) \\
& + 16 f(\textbf{x} + \eta e_i)
- 2 f(\textbf{x} + 2 \eta e_i) \; \big)
\end{split}
\end{equation}

\noindent In formula \ref{eq::fourth_order_diff}, the index $i=\left\lbrace 1,...,NR \right\rbrace$ is the index of the variables for which the derivative is to be evaluated. The variable $e_i$ is the $i-$th unit vector. The term $\eta$ is the perturbation and it is fixed small enough to achieve the convergence of the iterative process. 

Computing the exact inverse of the Hessian matrix arises numerical difficulties. However, as described by Wright et al. in \cite{wright1999numerical}, the Newton algorithm does not require a complete knowledge of the Hessian matrix. During the computation of the inverse of the Hessian matrix, the Hessian matrix is multiplied with a descent direction vector resulting in a vector. Therefore, only the results of the Hessian vector product is required. Using the Taylor expansion, this product is equal to the equation \ref{eq::hesvectdot}
\begin{equation} \label{eq::hesvectdot}
\nabla^2 f(\textbf{x}) \: \textbf{p} = \dfrac{\nabla f(\textbf{x}+\eta  \: \textbf{p})-\nabla f(\textbf{x})}{\eta}
\end{equation}
with $\eta$ the perturbation and $\textbf{p}$ the descent direction vector, fixed equal to the gradient at initialization. As a result, the extensive computation of the full Hessian matrix is bypassed using only the gradient. Finally, the complete ApHeN resolution scheme is presented in the algorithm \ref{algo:newcg}.\\

%
%
%
\textbf{Theoretical convergence rate}
APHEN is based on Newton's iterative method but it relies on an approximation of the Hessian matrix instead of the exact Hessian matrix. The reason is that although the exact Newton's method convergence is quadratic \cite{wright1999numerical}, the computation of the exact Hessian matrix is too time consuming for tensor application. Therefore, APHEN has a superlinear convergence such that 
\begin{equation} \label{eq::superlinconv}
\underset{n\rightarrow \infty}{lim} \dfrac{\left \| \textbf{B}_n - \nabla^2 f(\textbf{x}^*)\textbf{p}_n \right \|}{\left \| \textbf{p}_n \right \|} = 0
\end{equation}
with $\textbf{x}^*$ the point of convergence, $\textbf{p}_n$ the search direction and $\textbf{B}_n$ the approximation of the Hessian matrix. Practically, the convergence rate is described the equation below. 
\begin{equation} \label{eq::convergencerate}
q \approx 
\log \dfrac{|\textbf{x}_{n+1}-\textbf{x}_{n}|}{|\textbf{x}_{n}-\textbf{x}_{n-1}|} ( \log \dfrac{|\textbf{x}_{n}-\textbf{x}_{n-1}|}{|\textbf{x}_{n-1}-\textbf{x}_{n-2}|} )^{-1}
\end{equation}

\SetAlFnt{\footnotesize}
\SetAlCapFnt{\footnotesize}
\SetAlCapNameFnt{\footnotesize}

\begin{algorithm}[b!]
\setstretch{1.25}
\DontPrintSemicolon

\KwData{tensor $\mathscr{X} \in \mathbb{R}^{I\times J\times K}$, latent factors $(P,Q)$}

\KwResult{$\textbf{A}, \mathscr{D}^A, \textbf{H}, \mathscr{D}^B, \textbf{B}$ from tensor decomposition}

\Begin{


random initialization A$\in\mathbb{R}^{I\times P}$

random initialization H$\in\mathbb{R}^{P\times Q}$

random initialization B$\in\mathbb{R}^{J\times Q}$

set $\mathscr{D}^A_k \in\mathbb{R}^{P\times P}$ equal to 1 for $k=1,...,K$

set $\mathscr{D}^B_k \in\mathbb{R}^{Q\times Q}$ equal to 1 for $k=1,...,K$

\textbf{x}$_0$ $\gets$ flatten($\textbf{A}$, $\mathscr{D}^A$, $\textbf{H}$, $\mathscr{D}^B$, $\textbf{B}$) as described in (\ref{eq::vectorize_x})

$n=0$

\tcc{Error Minimization Loop}
\Repeat{maximum number of iterations or stopping criteria}
{
$\nabla f_n=\dfrac{\partial}{\partial x_i} f(\textbf{x}_n)$ $\longleftarrow$ gradient of $f$ at $\textbf{x}_n$ with (\ref{eq::fourth_order_diff})

$\textbf{p}_0 = - \nabla f_n \gets$ initial descent direction

$m = 0$

$\nabla f_{n_2}=\dfrac{\partial}{\partial x_i} f(\textbf{x}_n+\eta \textbf{p}_m)$

$\nabla^2 f(\textbf{x}_n)\textbf{p}_m = \dfrac{\nabla f_{n_2} - \nabla f_{n}}{\eta}$

\tcc{Search Direction CG Loop}
\Repeat{maximum number of iterations or stopping criteria}
{
\tcc{update rules as described in \cite{wright1999numerical}}
$\textbf{p}_m$ $\gets$ CG method applied to $\nabla^2 f(\textbf{x}_n)\textbf{p}_m = -\nabla f_n$ to determine the search direction $\textbf{p}_m$

m = m + 1
}

$\textbf{p}_n \gets \textbf{p}_m$


$\alpha_n=\underset{\textbf{x}_n, \nabla f_n, \textbf{p}_n}{\text{argmin}} f$ $ \gets$ Wolfe's line search for optimal step size

$\textbf{x}_{n+1} = \textbf{x}_n + \alpha_n \textbf{p}_n$
    
$n = n +1$

}

\KwRet{$\textbf{A}, \mathscr{D}^A, \textbf{H}, \mathscr{D}^B, \textbf{B}$}
}

\caption{ApHeN algorithm applied to Paratuck2 decomposition for a tensor $\mathscr{X}\in\mathbb{R}^{I\times J\times K}$ of latent factors $(P,Q)$\label{algo:newcg}}
\end{algorithm}
\subsection{Latent Predictions on Paratuck2 Tensor Decomposition }
Besides a Paratuck2 application of ApHeN for user-device authentication, our contribution resides in latent predictions. 
Following the tensor decomposition, latent variables are highlighted but modeled only past interactions. Hereinafter, the aim is to leverage  past information to predict the users' authentication. We briefly describe machine learning regression and neural networks used in our experiments. 

Decision Trees (DT) are a widely used machine learning technique \cite{lior2014data}. They are used to predict the value of a variable by learning simple decision rules from the data \cite{kim2015decision,breiman2017classification}. However, their regression decision rules have some limitations. Therefore, outpacing DT capabilities, neural networks including Multi-Layer Perceptron (MLP), Convolutional Neural Network (CNN) and Long-Short-Term-Memory (LSTM), and their applications, have skyrocketed for the past few years \cite{liu2017survey}. MLP consists of at least three layers: one input layer, one output layer and one or more hidden layer \cite{goodfellow2016deep}. Each neuron of the hidden layer transforms the values of the previous layer with a non-linear activation function. Although MLP is applied in deep learning, it lacks the possibility of modeling short term and long term events. This feature is found in LSTM \cite{goodfellow2016deep}. The LSTM has a memory block connected to the input gate and the output gate. The memory block is activated through a forget gate, resetting the memory information. However, for classification and computer vision, CNN is worth considering. In a CNN, the neurons are capable of extracting high order features in successive layers \cite{hubel1959receptive}. Through proper classification, the CNN is able to detect and predict various tasks including activities recognition \cite{abdel2012applying, zheng2014time}. 

\section{Experiments}
First, we highlight the numerical advantages of APHEN in comparison to other popular numerical schemes. Secondly, we rely on APHEN for the interactions modeling of the user-device authentication and their predictions.

\subsection{APHEN vs Other Numerical Schemes}
Hereinafter, we investigate the convergence behavior of APHEN in comparison to other numerical resolution methods. First, we define the concept of convergence rate and convergence speed. Then, we compare APHEN with 6 different algorithms applied to Paratuck2: 
\begin{itemize}
\item ALS, Alternating Least Squares \cite{bro1998multi,charlier2018non}
\item GD, Gradient Descent \cite{wright1999numerical}
\item NAG, Nesterov Accelerated Gradient \cite{nesterov2007gradient}
\item SAGA \cite{defazio2014saga}
\item Adam \cite{kingma2014adam}
\item BFGS \cite{broyden1970convergence,fletcher1970new,goldfarb1970family,shanno1970conditioning}
\end{itemize}
The simulations are conducted on a PC with an Intel Core i7 CPU and 16GB of RAM. All the resolution schemes have been implemented in Julia.\\

%
%
%
\textbf{Convergence speed definition}
The definition of the convergence rate in \ref{eq::convergencerate} does not illustrate the time evolution between each iteration. Therefore, we define two notions: the iteration-based convergence speed and the time-based convergence speed. The convergence speed is defined as the absolute value of the linear slope, denoted by $|a|$, of each curve shown afterward according to the equation $y=ax + b$. The bigger the convergence speed, the faster convergence and the lower the numerical errors in the tensor decomposition. Two types of convergence speed are characterized: the iteration-based convergence speed and the time-based convergence speed. The iteration-based convergence speed, and the time-based convergence speed, measures the evolution of the numerical errors according to the iterations, and to the time, respectively.\\

%
%
%
\textbf{Numerical Convergence Highlights}
Seven tensor sizes have been defined, 5$\times$5$\times$5, 10$\times$10$\times$10, 15$\times$10$\times$10, 15$\times$15$\times$15, 25$\times$20$\times$15, 50$\times$40$\times$20 and 100$\times$100$\times$20, with respective latent factors $(2,3)$, $(3,4)$, $(5,4)$, $(5,6)$, $(10,9)$, $(15,14)$ and $(3,5)$. The tensor dimensions and the latent factors have been chosen arbitrarily since the experiments have shown similar results for any tensor for any combination of latent factors. Each tensor entry is incremented by one in comparison to the previous entry, with the initial entry fixed to one. Additionally, for all the simulation, we define the convergence criterion such that $|f(\textbf{x}_n)-f(\textbf{x}_{n-1})|\left(|f(\textbf{x}_n)|\right)^{-1}< 10^{-6}$.

\begin{figure}[h]
  \centering
  \includegraphics[scale=0.4]{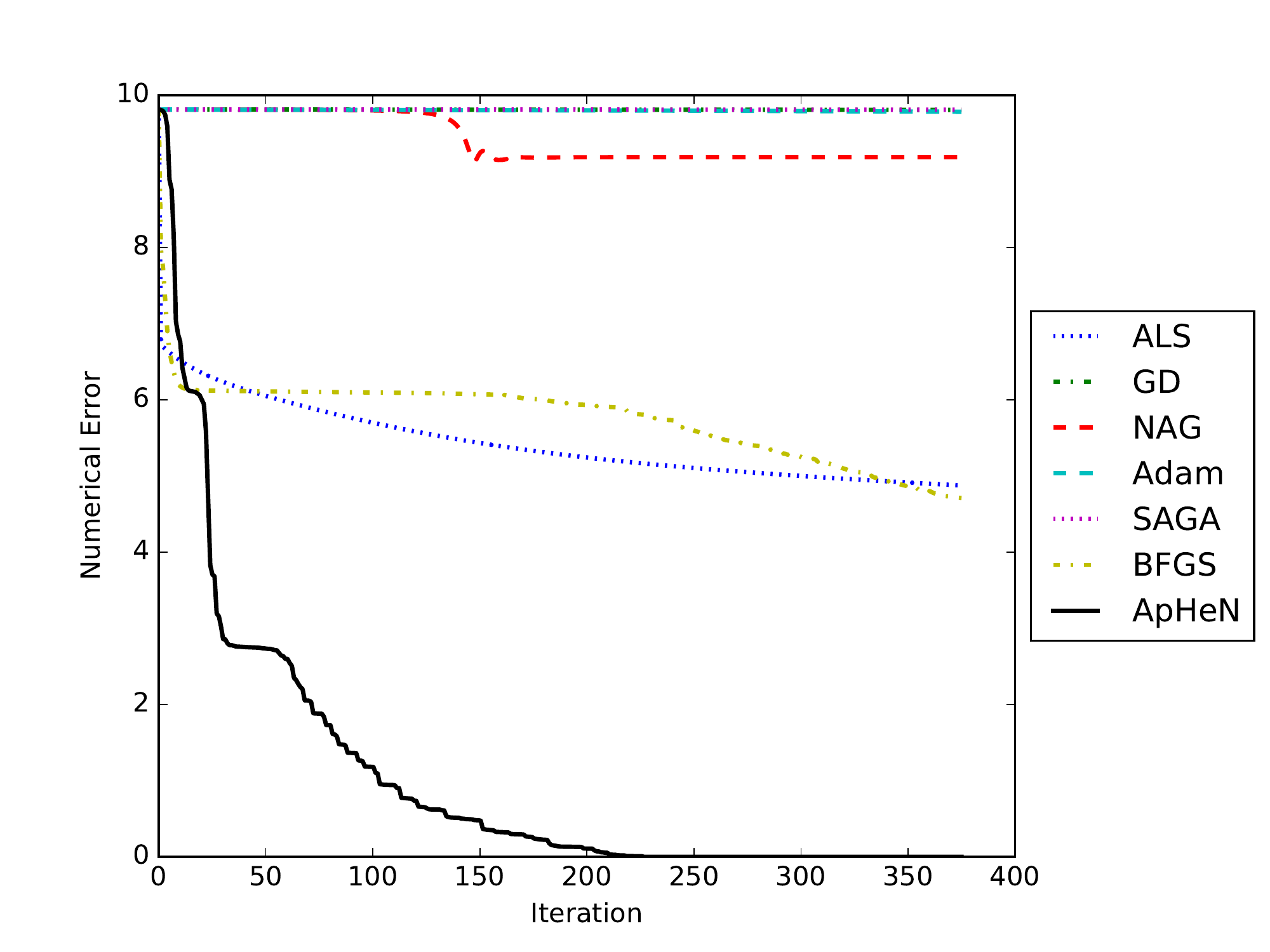}
  \includegraphics[scale=0.4]{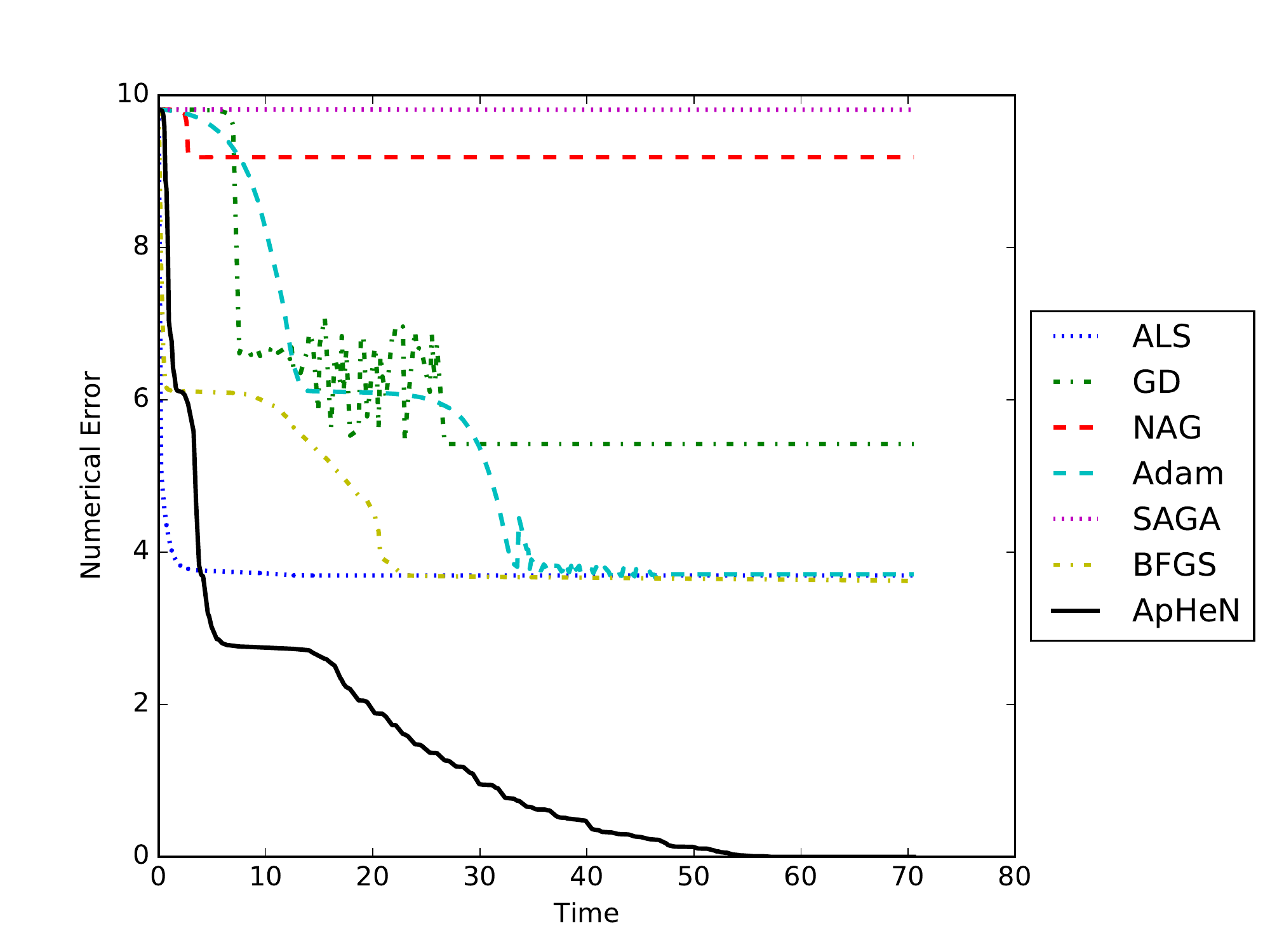}
  \caption{Iteration-based and time-based convergences of the different numerical resolution schemes applied to Paratuck2 for a tensor of size 10$\times$10$\times$10 with latent factors $(3,4)$}
  \label{fig::twodplot}
\end{figure}

Figure \ref{fig::twodplot} highlights the convergence speeds of the different resolution schemes for a tensor of size 10$\times$10$\times$10 with latent factors $(3,4)$. ALS, BFGS and APHEN show significant superior convergence speeds, both iteration-based and time-based. The ALS scheme decreases the fastest at the beginning of the process but it fails rapidly to determine the solution having the lowest numerical errors. Although, APHEN has slightly longer computation time than ALS, it is the only method capable of determining the optimal solution. Surprisingly, all the gradient schemes have significantly lower convergence speeds including Adam. 

\begin{figure}[h]
  \centering
  \includegraphics[scale=0.4]{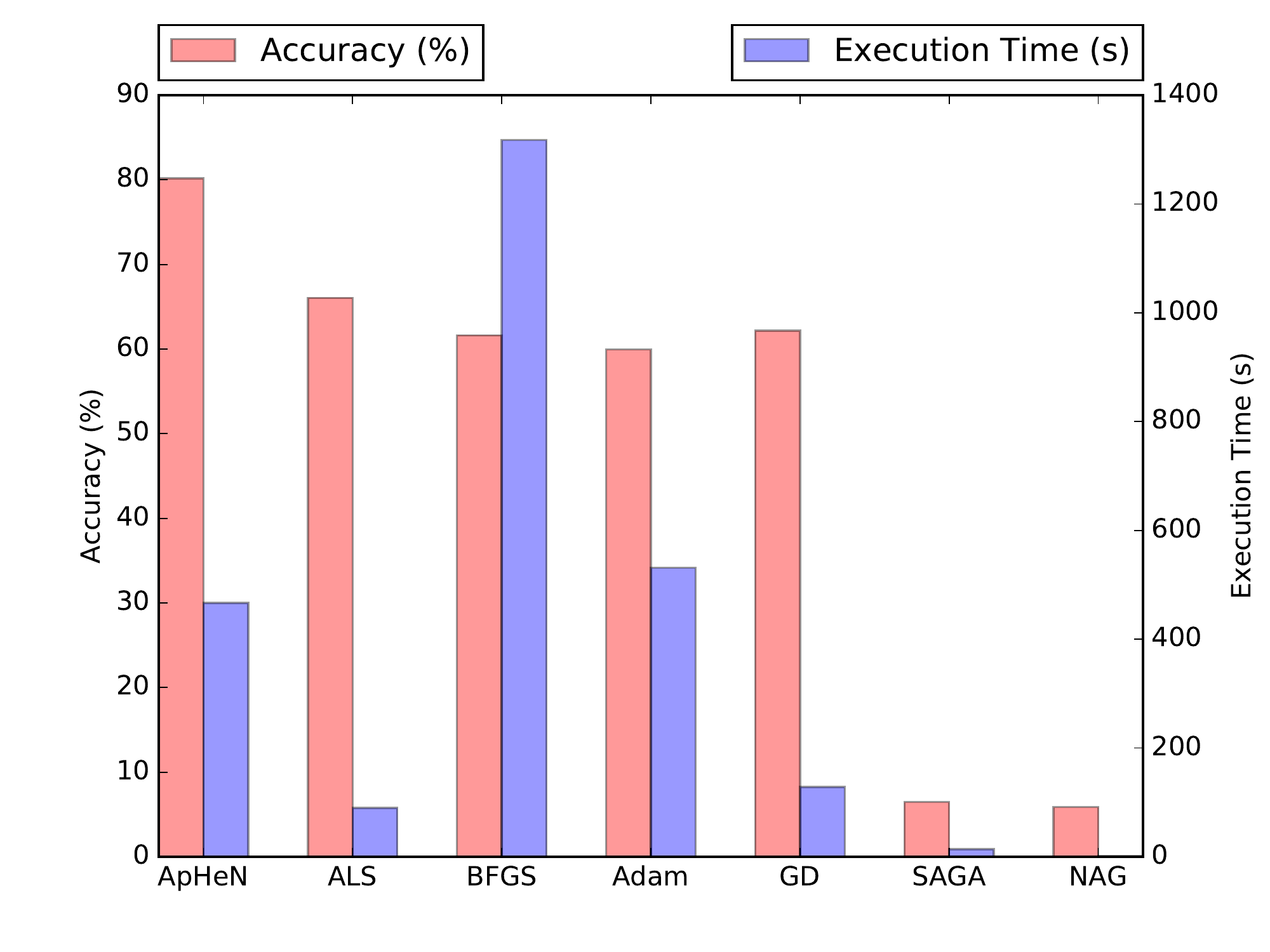}
  \caption{Accuracy of each resolution method (left column) with their respective execution time (right column) at convergence applied to Paratuck2 for a tensor of size 15$\times$15$\times$15 with latent factors $(5,6)$}
  \label{fig::barplot}
\end{figure}

Figure \ref{fig::barplot} highlights the accuracy and the execution time of the different resolution schemes for a tensor of size 15$\times$15$\times$15. The accuracy is defined as 
$100(1 - \frac{\log ||\mathscr{X}-\mathscr{\hat{X}}||}{\log ||\mathscr{X}||}1_{||\mathscr{X}-\mathscr{\hat{X}}||>1} )$. APHEN has the best accuracy followed by ALS, BFGS, Adam and the other schemes. BFGS and Adam have longer time of execution than APHEN for a significantly lower accuracy at convergence. ALS has a slightly faster execution time than APHEN but APHEN has a better accuracy at convergence. 

These graphical results are completed by the tables \ref{tab::convergencehighlights} and \ref{tab::accuracy}. The table \ref{tab::convergencehighlights} shows APHEN has the fastest convergence speeds in all the simulation followed closely by ALS, BFGS and Adam. The performance of the other schemes are significantly lower. Furthermore, the table \ref{tab::accuracy} highlights the superiority of APHEN to determine the solution having the lowest numerical residual errors at the convergence of the calculation.

To summarize, we showed APHEN provides faster convergence speeds and lower residual errors for similar execution time. Thus, we use APHEN with Paratuck2 in our work.

\begin{table*}[t]
\centering
\caption{Convergence speed highlights of the different resolution schemes (bigger is better).}
\label{tab::convergencehighlights}
\begin{tabular}{cccccccccc}
\toprule
Convergence & Tensor & Latent & Grad.-Free & \multicolumn{4}{c}{Hessian-Free} & \multicolumn{2}{c}{Hessian Approximation} \\
Type & Size & Factors & ALS & GD & NAG & Adam & SAGA & BFGS & APHEN\\
\midrule
Iteration & 5$\times$5$\times$5 & 2$\times$3 & 0.0367 & 0.0001 & 0.0059 & 0.0002 & 0.0002 & 0.0322 & \textbf{0.0676}\\
Iteration & 10$\times$10$\times$10 & 3$\times$4 & 0.0231 & 0.0001 & 0.0028 & 0.0001 & 0.0042 & 0.0196 & \textbf{0.0490}\\
Iteration & 15$\times$10$\times$10 & 5$\times$4 & 0.0212 & 0.0001 & 0.0029 & 0.0001 & 0.0043 & 0.0238 & \textbf{0.0451}\\
Iteration & 15$\times$15$\times$15 & 5$\times$6 & 0.0135 & 0.0096 & 0.0001 & 0.0001 & 0.0001 & 0.0136 & \textbf{0.0268}\\
Iteration & 25$\times$20$\times$15 & 10$\times$9 & 0.0136 & 0.0045 & 0.0001 & 0.0109 & 0.0001 & 0.0146 & \textbf{0.0250}\\
Iteration & 50$\times$40$\times$20 & 15$\times$14 & 0.0132 & 0.0001 & 0.0017 & 0.0001 & 0.0001 & 0.0137 & \textbf{0.0232}\\
Iteration & 100$\times$100$\times$20 & 3$\times$5 & 0.0856 & 0.0001 & 0.0001 & 0.0001 & 0.0001 & 0.0913 & \textbf{0.1032}\\

Time & 5$\times$5$\times$5 & 2$\times$3 & 0.2706 & 0.1543 & 0.029 & 0.1514 & 0.0335 & 0.2364 & \textbf{0.2883}\\
Time & 10$\times$10$\times$10 & 3$\times$4 & 0.0245 & 0.0157 & 0.0025 & 0.0152 & 0.0033 & 0.0245 & \textbf{0.0339} \\
Time & 15$\times$10$\times$10 & 5$\times$4 & 0.0146 & 0.0107 & 0.0001 & 0.0076 & 0.0017 & 0.013 & \textbf{0.0156}\\
Time & 15$\times$15$\times$15 & 5$\times$6 & 0.0055 & 0.0017 & 0.0001 & 0.0001 & 0.0005 & 0.0049 & \textbf{0.0058} \\
Time & 25$\times$20$\times$15 & 10$\times$9 & 0.0042 & 0.002 & 0.0001 & 0.0020 & 0.0001 & 0.0037 & \textbf{0.0044}\\
Time & 50$\times$40$\times$20 & 15$\times$14 & 0.0020 & 0.0015 & 0.0002 & 0.0020 & 0.0002 & 0.0022 &\textbf{0.0033}\\
Time & 100$\times$100$\times$20 & 3$\times$5 & 0.0013 & 0.0005 &  0.0001 & 0.0001 & 0.0001 & 0.0013 & \textbf{0.0015}\\
\bottomrule
\end{tabular}
\end{table*}


\begin{table*}[t]
\caption{Accuracy of the different resolution schemes at convergence (bigger is better).}
\centering
\label{tab::accuracy}
\begin{tabular}{ccccccccc}
\toprule
Tensor & Latent & Grad.-Free & \multicolumn{4}{c}{Hessian-Free} & \multicolumn{2}{c}{Hessian Approximation} \\
Size & Factors & ALS & GD & NAG & Adam & SAGA & BFGS & APHEN\\
\midrule
5$\times$5$\times$5 & 2$\times$3 & 80.8952 & 73.3299 & 8.6596 & 75.0865 & 10.0273 & 79.2534 & \textbf{99.9999}\\
10$\times$10$\times$10 & 3$\times$4 & 62.3542 & 50.3829 & 6.9655 & 62.7238 & 10.2185 & 66.3565 & \textbf{98.8941}\\
15$\times$10$\times$10 & 5$\times$4 & 70.2395 & 64.4161 & 5.8585 & 65.3552 & 11.2226 & 66.0316 & \textbf{88.0582}\\
15$\times$15$\times$15 & 5$\times$6 & 66.0539 & 62.1681 & 5.8687 & 59.9574 & 6.4757 & 61.6300 & \textbf{80.1714}\\
25$\times$20$\times$15 & 10$\times$9 & 65.4691 & 43.6754 & 4.5421 & 44.1859 & 10.5425 & 57.0254 & \textbf{68.5696} \\
50$\times$40$\times$20 & 15$\times$14 & 72.3543 & 50.3830 & 6.3749 & 62.7238 & 7.2364 & 66.3566 & \textbf{87.9709} \\
100$\times$100$\times$20 & 3$\times$5 & 49.4730 & 38.7512 & 2.8462 & 48.1267 & 3.6195 & 49.3348 & \textbf{55.7678} \\
\bottomrule
\end{tabular}
\end{table*}

\subsection{User-Device Authentication Monitoring for Financial Recommendation}
First, we discuss the completion and the resolution of the tensor. Secondly, we rely on the results of Paratuck2 for the predictions of the authentication with the neural networks.\\

%
%
\textbf{User-Computer Authentication and Data Availability}
For the sake of the reproducibility of the experiments, we present the approach with a public data set. In 2014, the Los Alamos National Laboratory enterprise network published the anonymized user-computer authentication logs of their laboratory \cite{hagberg-2014-credential}, and available at {\color{blue} \underline{\url{https://csr.lanl.gov/data/auth/}}}. Each authentication event is composed of the authentication time (in Unix time), the computer label and the user label such as, for instance, "1,U1,C1". In total, more than 11,000 users and 22,000 computers are listed representing 13 GB of data. \\

\textbf{Construction of the user-computer authentication tensor}
We randomly select 150 users and 300 computers within the dataset representing more than 60 millions lines. The first two months of authentication events have been compressed into 50 time intervals, corresponding to 25 working days per month. A tensor $\mathscr{X}\in \mathbb{R}^{I\times J\times K}$ of size of 150$\times$300$\times$50 is built. The first dimension, denoted by $I$, represents the users, the second dimension, denoted by $J$, the computers and the last dimension, $K$, stands for the time intervals. \\


\textbf{Limitations of the CP decomposition}
The CP decomposition expresses the original tensor into a sum of rank one tensors. Therefore, the user-computer authentication tensor is decomposed as a sum of user-computer-time rank-one tensors.
However, in the case of strong imbalance, CP leads to underfitting or overfitting one of the dimension \cite{acar2011scalable}. Within the dataset, we can find 2 users that connect to at least 20 different computers. Therefore, a rank equal to 2, one per user, underfits the computer connections. A rank equal to 20, one per machine, overfits the number of users. In the table \ref{tab::overfitcp}, the underfitting is underlined by significant residual errors at convergence. The overfitting is detected by a good understanding of the data since the residual errors tend to be small. Hence, the Paratuck2 decomposition is chosen to model properly each dimension of the original tensor. \\


\begin{table}[t]
  \caption{In CP, for imbalanced dataset, underfitting one dimension is highlighted by significant residual errors. Overfitting is difficult to measure because of the low residual errors. A good understanding of the data is required to estimate it.}
  \label{tab::overfitcp}
  \centering
    \begin{tabular}{cccc}
	\toprule
    Tensor Size & Rank & Residual Errors & $\frac{|f(x_n) - f(x_{n-1})|}{|f(x_n)|} $\\
    \midrule
    2$\times$20$\times$30 & 2 & \textbf{50.275} & $<10^{-6}$ \\
    2$\times$20$\times$30 & 20 & \textbf{1.147} & $<10^{-6}$ \\
    \bottomrule
  \end{tabular}
\end{table}

\begin{figure}[t]
\begin{center}
\begin{tikzpicture}[scale=0.625, transform shape]
\draw  (-4.5,2) rectangle (-3,-1);
\fill [gray, opacity=.4]  (-4.5,2) rectangle (-3,-1);
\draw (-3.75,0.75) node[below]{\large{$A$}};

\draw (2.5,-0.5) rectangle (4.5,1.5) node (v1) {};
\draw (4.5,1.5) -- (5,2);
\draw (2.5,1.5) -- (3,2);
\draw (4.5,-0.5) -- (5,0);
\draw (3,2) -- (5,2);
\draw (5,0) -- (5,2);
\draw [dashdotted] (3,0) -- (5,0);
\draw [dashdotted] (3,0) -- (3,2);
\draw [dashdotted] (2.5,-0.5) -- (3,0);
\fill [gray, opacity=.4] (4.5,-0.5) -- (5,0) -- (3,2) -- (2.5,1.5)  -- cycle;
\draw (-1.5,1.3) node[below]{\large{$D^A$}};

\draw (0,2) rectangle (2,0.5);
\fill [gray, opacity=.4] (0,2) rectangle (2,0.5);
\draw (1,1.4) node[below]{\large{$H$}};

\draw (-2.55,-0.05) rectangle (-1,1.5) node (v1) {};
\draw (-1,1.5) -- (-0.5,2);
\draw (-2.55,1.5) -- (-2.05,2);
\draw (-1,-0.05) -- (-0.5,0.45);
\draw (-2.05,2) -- (-0.5,2);
\draw (-0.5,0.45) -- (-0.5,2);
\draw [dashdotted] (-2.05,0.45) -- (-0.5,0.45);
\draw [dashdotted] (-2.05,0.45) -- (-2.05,2);
\draw [dashdotted] (-2.55,-0.05) -- (-2.05,0.45);
\fill [gray, opacity=.4] (-1,-0.05) -- (-0.5,0.45) -- (-2.05,2.0) -- (-2.55,1.5)  -- cycle;
\draw (3.8,1.05) node[below]{\large{$D^B$}};

\draw  (5.5,2) rectangle (8,0);
\fill [gray, opacity=.4] (5.5,2) rectangle (8,0);
\draw (6.8,1.25) node[below]{\large{$B^T$}};

\draw[decoration={brace,raise=5pt, amplitude=7pt},decorate] (-4.5,2.0) -- node[below=-35pt] {\textit{into P groups}} (-3,2.0);
\node at (-3.75,3.25) {\textit{clustering I users}};

\draw[decoration={brace,mirror,raise=5pt, amplitude=10pt},decorate] (-2.5,-0.10) -- (-0.5,-0.10);
\node at (-1.5,-1.3) {\textit{Evolution over time of}};
\node at (-1.5,-1.65) {\textit{P users groups connections}};

\draw[decoration={brace,raise=5pt, amplitude=7pt},decorate] (0,2.0) -- node[below=-35pt] {\textit{between P and Q latent comp.}} (2,2.0);
\node at (1,3.25) {\textit{static asymmetry}};

\draw[decoration={brace,mirror,raise=5pt, amplitude=7pt},decorate] (2.5,-0.40) -- (5,-0.40);
\node at (3.75,-1.3) {\textit{Evolution over time of }};
\node at (3.75,-1.65) {\textit{Q computers groups connections}};

\draw[decoration={brace,raise=5pt, amplitude=7pt},decorate] (5.5,2.0) -- node[below=-35pt] {\textit{into Q groups}} (8,2.0);
\node at (6.75,3.25) {\textit{clustering J computers}};

\end{tikzpicture}
\caption{Paratuck2 decomposition applied to user-computer authentication. The neural network predictions are performed on the tensor $\mathscr{D}^A$.}
\label{fig::PARATUCK2_exp}
\end{center}
\end{figure}
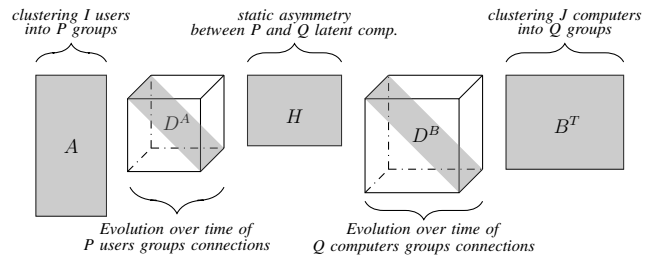

\textbf{Paratuck2 Tensor Resolution}
Paratuck2 decomposes the main tensor $\mathscr{X}\in\mathbb{R}^{I\times J \times K}$ into a product of matrices and sparse tensors as shown in the figure \ref{fig::PARATUCK2_exp}. The matrix \textbf{A} factorizes the users into $P$ groups. We observe 15 different groups of users, and therefore, $P$ equals to 15. The sparse tensor $\mathscr{D}^A$ reflects the temporal evolution of the connections of the $P$ users groups. The matrix $\textbf{H}$ represents the asymmetry between the $P$ users groups and the $Q$ computers groups. We notice 25 different groups of machines related to different authentication profiles, and consequently, $Q$ equals to 25. The sparse tensor $\mathscr{D}^B$ illustrates the temporal evolution of the connections of the $Q$ computers groups. Finally, the matrix $\textbf{B}$ factorizes the computers into $Q$ latent groups of computers. \\


\textbf{Latent Predictions for Financial Recommendation}
To achieve higher subscription rates during the advertising campaign of financial products, we explore the latent predictions for targeted recommendation based on the future user-computer authentication. The results of Paratuck2 contain the users' temporal information and the computers' temporal information in the sparse tensors $\mathscr{D}^A$ and $\mathscr{D}^B$, respectively. Predicting the users' authentication allows the banks to build a more complete financial awareness profile of their clients for optimized advertisement. 

\begin{figure}[h]
  \centering
 \includegraphics[scale=0.4]{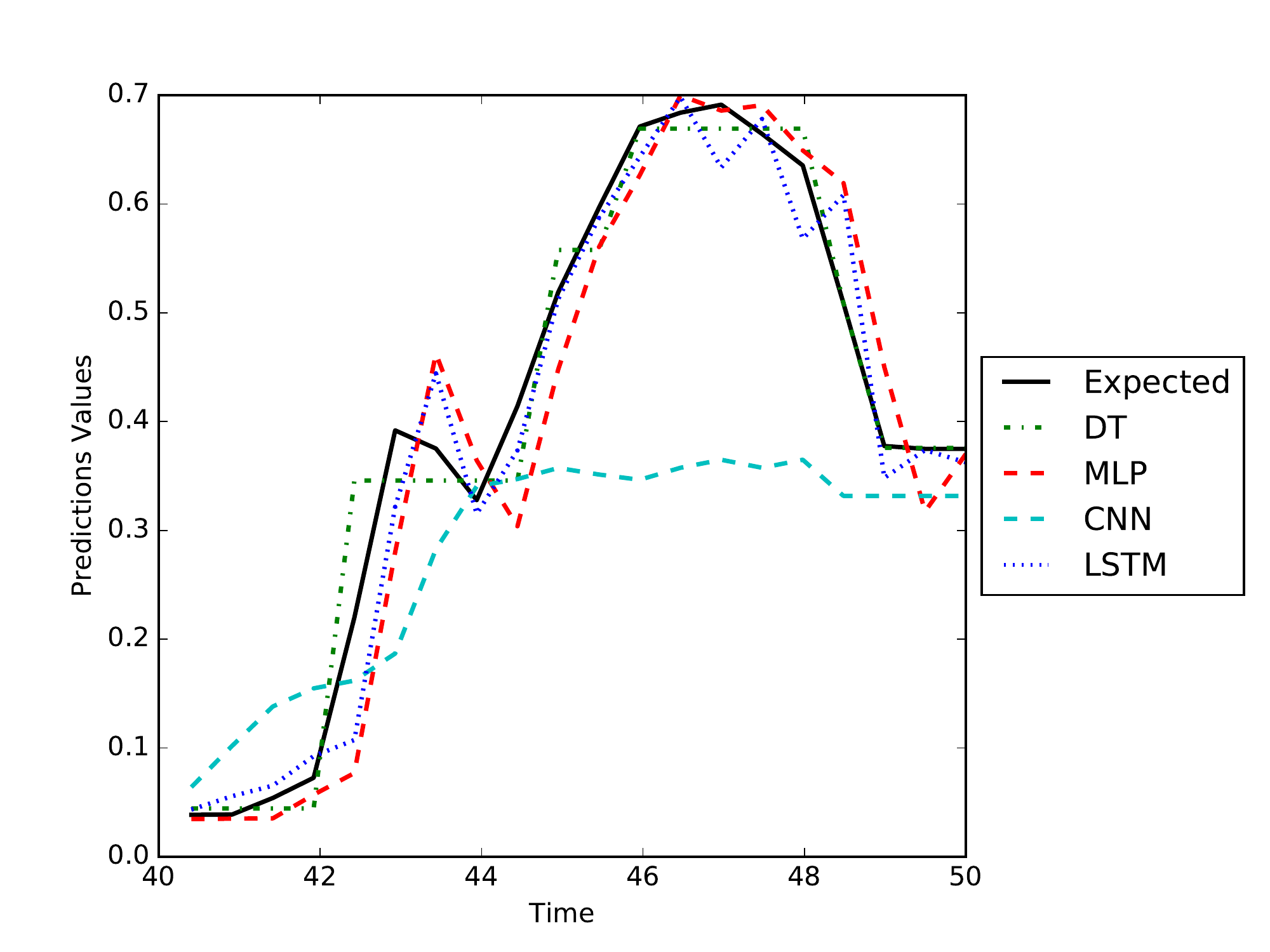}
  \caption{Two weeks prediction of the evolution of the latent users' authentication according to the different models used}
  \label{fig::preds}
\end{figure}

In figure \ref{fig::preds}, we highlight the results of the predictions of the users' authentication for a specific group of clients, corresponding with one specific latent factor $P$. Four different methods have been used for the predictions, DT, MLP, CNN and LSTM. All methods have been trained on a six weeks period. Then, the users' authentication for the next two weeks are predicted with a rolling time window of one day. The figure \ref{fig::preds} highlights visually that the LSTM models the most accurately the future users' authentication. It is followed by the MLP, the DT, and finally the CNN. We underline this preliminary statement using six well-known error measures. The Mean Absolute Error (MAE), the Mean Directional Accuracy (MDA), the Pearson correlation, the Jaccard distance, the cosine similarity and the Root Mean Square Error (RMSE) are used to determine objectively the most accurate predictive method. 
The table \ref{tab::predictionerrors1} describes the error measures related to the figure \ref{fig::preds}. As previously seen, the LSTM is the closest to the true authentication since it has the lowest error values. Then, the MLP comes second, the DT third, and the CNN last.  

To conclude, with the aim to better target the clients that might be interested by financial products during the bank's advertising campaigns, we can conclude that LSTM combined with Paratuck2 models the best the future users' authentication. As the majority of the user's authentication are sequence-based, it is legitimate to find out LSTM gives the best results for the predictions. Effectively, each user has a recurrent pattern in the authentication process depending on its activities of the day. Therefore, by using APHEN for Paratuck2 and LSTM for predictions, the bank gain a very competitive advantage for the personalized products recommendation, based only on its clients' authentication on the mobile application.

\begin{table}[h]
\caption{Latent predictions errors on the users' authentication with decision tree and neural networks}
\centering
\label{tab::predictionerrors1}
\begin{tabular}{ccccc}
\toprule
Error Measure & DT & MLP & CNN & LSTM \\
\midrule
MAE & 0.0965 & 0.0506 & 0.1106 & \textbf{0.0379} \\
MDA & 0.1579 & 0.7447 & 0.5263 & \textbf{0.6842} \\
Pearson corr. & 0.8537 & 0.9598 & 0.8885 & \textbf{0.9753} \\
Jaccard dist. & 0.2257 & 0.1206 & 0.2648 & \textbf{0.0911} \\
cosine sim. & 0.9587 & 0.9891 & 0.9745 & \textbf{0.9914} \\
RMSE & 0.1306 & 0.0695 & 0.3140 & \textbf{0.0477} \\
\bottomrule
\end{tabular}
\end{table}

\section{Conclusion And Future Work}
In this paper, we presented an Hessian-based algorithm, APHEN, that does not require a full knowledge of the Hessian matrix. It was applied to resolution of the Paratuck2 tensor decomposition. APHEN reduces at a minimum the numerical errors inherited from the tensor decomposition. Furthermore, it has higher convergence speeds than other popular methods such as NAG or Adam. We used derivatives approximation evaluated with finite difference schemes to propose an accessible framework for all tensor decompositions. The experiments were conducted on tensors of different sizes with different latent factors. Additionally, we showcased an application in the context of mobile banking application. We used Paratuck2 and state of the art machine learning and neural networks to profile and predict the latent users' authentication. By modeling the clients' past and future authentication on their mobile application, the banks are able to build a financial awareness profile of their clients to advert different types of products. The banks have realized the promising potential of the clients' digital behavior to face the increasing competition coming from the new regulation directives.

As future work, we plan on showing the versatility of APHEN to all tensor decompositions. We will compare APHEN's performance  for all existing tensor decomposition against the other existing tensor resolution algorithms specific to each tensor decomposition. Then, we will assess the influence of the line search and the performance of adaptive line searches while improving the GPU compatibility of the algorithm to increase the size of the experiments. Finally, the financial recommendation depending of the user-device authentication on a mobile banking application will be further extended. The navigation usage, the time gap between each action and the type of device used will be monitored to further improve the bank's advertising campaigns of their products to the appropriate clients. 

\bibliographystyle{./IEEEtran}
\bibliography{./IEEEabrv,./IEEEexample}

\end{document}